\documentstyle [12pt]{article}

\title {A model theoretic Baire category theorem for simple theories and its applications}
\author {Ziv Shami}
\newtheorem {theorem}{Theorem}[section]
\newtheorem {lemma}[theorem]{Lemma}
\newtheorem {definition}[theorem]{Definition}

\newtheorem {fact}[theorem]{Fact}

\newtheorem {corollary}[theorem]{Corollary}

\newtheorem {remark}[theorem]{Remark}

\newtheorem {subclaim}[theorem]{Subclaim}

\def\proof {\noindent \textbf{Proof:} }

\newsavebox{\indbin}
\savebox{\indbin}{\begin{picture}(0,0)
\newlength{\gnu}
\settowidth{\gnu}{$\smile$} \setlength{\unitlength}{.5\gnu} \put(-1,-.65){$\smile$}
\put(-.25,.1){$|$}
\end{picture}}
\newcommand{\nonfork}[3]
{\mbox{$\begin{array}{ccc} \mbox{$#1$} & \usebox{\indbin} & \mbox{$#2$} \\
        & \mbox{$#3$} &
\end{array}$}}

\newsavebox{\sindbin}
\savebox{\sindbin}{\begin{picture}(0,0)
\newlength{\sgnu}
\settowidth{\sgnu}{$\smile$} \setlength{\unitlength}{.5\sgnu} \put(-1,-.65){$\smile$}
\put(-.25,.1){$|s$}
\end{picture}}

\newcommand{\snonforkempty}[2]
{\mbox{$\begin{array}{ccc} \mbox{$#1$} & \usebox{\sindbin} & \mbox{$#2$}
\end{array}$}}
\newcommand{\sfork}[3]
{\mbox{$\begin{array}{ccc} \mbox{$#1$} & \!\mbox{$\!\!\not\!\:\usebox{\sindbin}$} & \mbox{$#2$} \\
        & \mbox{$#3$} &
\end{array}$}}
\newcommand{\sforkempty}[2]
{\mbox{$\begin{array}{ccc} \mbox{$#1$} & \!\mbox{$\!\!\not\!\:\usebox{\sindbin}$} & \mbox{$#2$}
\end{array}$}}

\def\card #1 {{\vert #1 \vert}}

\def\CC {{\cal C}}
\def\RR {{\cal R}}

\def\UU {{\cal U}}

\def\FF {{\cal F}}
\def\SS {{\cal S}}

\begin{document}
\maketitle

\begin{abstract}
We prove a model theoretic Baire category theorem for $\tilde\tau_{low}^f$-sets in a countable
simple theory in which the extension property is first-order and show some of its applications. We
also prove a trichotomy for minimal types in countable nfcp theories: either every type that is
internal in a minimal type is essentially-1-based by means of the forking topology  or $T$
interprets an infinite definable 1-based group of finite $D$-rank or $T$ interprets a
strongly-minimal formula.
\end{abstract}

\section{Introduction}
The goal of this paper is to generalize a result from [S1] and to give some applications. In [S1]
the first step for proving supersimplicity of countable unidimensional simple theories eliminating
hyperimaginaries is to show the existence of an unbounded type-definable forking-open set (a set
defined in terms of forking by formulas, see Definition \ref{tau definition}) of bounded finite
$SU_{se}$-rank (for definition see Section 4). In this paper we develop a general framework for
this kind of result. It is a new idea of a model theoretic Baire category theorem, namely, one
deals with certain "uniformly-definable" family of generalized closed sets (in complicated
"logic"); roughly speaking, given a partition of a complicated open set into countably many sets,
each of which is the intersection of a "uniformly definable" family of generalized closed sets, one
can find a forking-open set that is contained in some generalized closed set in one of these
families. So, the main point is that we obtain a very nice set (forking-open), but on the other
hand we can only require that it will be a subset of some generalized closed set in one of these
families and not in its intersection. In particular, it is not just the usual Baire category
theorem for a complicated topological space. The proof is quite similar to the proof in [S1] and
has some important consequences, e.g. in a countable wnfcp theory if for every non-algebraic
element $a$ (even in some fixed non-empty $\tilde\tau_{low}^f$-set) there is $a'\in
acl(a)\backslash acl(\emptyset)$ of finite $SU$-rank, then there exists a weakly-minimal formula.
We also prove a trichotomy for countable nfcp theories as indicated in the abstract.

We assume basic knowledge of simple theories. A good textbook on simple theories is [W]. The
notations follow usual conventions. $T$ will denote a complete first-order theory with no finite
models in some language $L$. We will work in some large saturated model $\cal C$ of $T$ (not
necessarily with elimination of imaginaries, unless stated otherwise). Ordinals will be denoted by
$\alpha$, $\beta$, $\gamma$,... Sets $A$, $B$, $C$,... will be small subsets of $\cal C$, i.e. of
cardinality strictly less than the cardinality of $\cal C$. The letters $a$, $b$, $c$,... denote
finite tuples from $\CC$ unless stated otherwise. $x$,$y$,$z$,... denote finite tuples of variables
unless stated otherwise. We use $p$, $q$, $r$,... to denote types (possibly partial) over some set.
For an invariant set $V$ (over some small set) and $n$, we denote by $V^n$ the set of $n$-tuples of
realizations of $V$.

\section{Preliminaries}
The forking topology is introduced in [S0] and is a variant of Hrushovski's and Pillay's topologies
from [H0] and [P], respectively. In this section $T$ is assumed to be simple and we work in a large
saturated model $\CC$ of $T$.

\begin{definition} \label{tau definition}\em
Let $A\subseteq \CC$ and let $x$ be a finite tuple of variables.\\
1) An invariant set $\UU$ over $A$ is said to be \em a basic $\tau^f$-open set over $A$ \em if
there is $\phi(x,y)\in L(A)$ such that $$\UU=\{a \vert \phi(a,y)\ \mbox{forks\ over}\ A \}.$$ Note
that the family of basic $\tau^f$-open sets over $A$ is closed under finite intersections, thus
forms a basis for a unique topology on $S_x(A)$. An open set in this topology is called a
$\tau^f$-open set over $A$ or a forking-open set over $A$.

\noindent 2) An invariant set $\UU$ over $A$ is said to be \em a basic $\tau^f_\infty$-open set
over $A$ \em if $\UU$ is a type-definable $\tau^f$-open set over $A$. The family of basic
$\tau^f_\infty$-open sets over $A$ is a basis for a unique topology on $S_x(A)$. An open set in
this topology is called a $\tau^f_\infty$-open set over $A$.
\end{definition}

Recall that a formula $\phi(x,y)\in L$ is \em low in $x$ \em if there exists $k<\omega$ such that
for every $\emptyset$-indiscernible sequence $(b_i \vert i<\omega)$, the set $\{\phi(x,b_i) \vert
i<\omega\}$ is inconsistent iff every subset of it of size $k$ is inconsistent. $T$ is low if every
$\phi(x,y)$ is low in $x$.

\begin{remark}\label{low_disjunction}\em
Assume $\phi(x,t)\in L$ is low in $t$ and $\psi(y,v)\in L$ is low in $v$ ($x\cap y$, $t\cap v$ may
not be $\emptyset$). Then $\theta(xy,tv)\equiv\phi(x,t)\vee \psi(y,v)$ is low in $tv$.
\end{remark}

\proof Let $k_1<\omega$ be a witness that $\phi(x,t)$ is low in $t$ and let $k_2<\omega$ be a
witness that $\psi(y,v)$ is low in $v$. Let $k=k_1+k_2-1$. By adding dummy variables we may assume
$x=y$ and $t=v$ (as tuples of variables). Let $(a_i \vert\ i<\omega)$ be indiscernible such that
$\{ \phi(a_i,t)\vee \psi(a_i,t)\vert i<\omega\}$ is inconsistent. Thus, every subset of $\{
\phi(a_i,t)\vert i<\omega\}$ of size $k_1$ is inconsistent, and every subset of $\{\psi(a_i,t)\vert
i<\omega\}$ of size $k_2$ is inconsistent. Thus every subset of size $k$ of
$\{ \phi(a_i,t)\vee \psi(a_i,t)\vert i<\omega\}$ is inconsistent.\\

Here we state some basic facts about the $\tau^f$-topology (See [S0, Claim 2.5], [S1, Remark 7.6]).

\begin{remark}\label{exists generic}\em
\noindent 1) The $\tau^f$-topology on $S_x(A)$ refines the Stone-topology of $S_x(A)$ for all
$x,A$.

\noindent 2) A basic  $\tau^f$-open set in a low theory is type-definable and every Stone-closed
subset of $(S_x(A),\tau^f)$ is a Baire topological space (i.e. the intersection of countably many
dense open sets in it is dense).

\noindent 3) Let $A$ be a small set. Let $F(x,y)$ be a type-definable relation over $A$ and let
$f(x)$ be an $A$-definable function. Let $\Gamma_{F,f}(x)=\exists y (F(x,y)\wedge
\nonfork{y}{f(x)}{A})$. Then $\Gamma_{F,f}(x)$ is $\tau^f$-closed over $A$ ([S0, Claim 2.5] is
slightly different, but the proof is the same).

\end{remark}

Recall the following definition from [S0] whose roots are in [H0].

\begin{definition}\label{projection closed}\em
We say that \em the $\tau^f$-topologies over $A$ are closed under projections ($T$ is PCFT over
$A$) \em if for every $\tau^f$-open set $\UU(x,y)$ over $A$ the set $\exists y \UU(x,y)$ is a
$\tau^f$-open set over $A$. We say that \em the $\tau^f$-topologies are closed under projections
($T$ is PCFT) \em if they are such over every set $A$.
\end{definition}

In [BPV, Proposition\ 4.5] the authors proved the following equivalence which, for convenience, we
will use as a definition (their definition involves extension with respect to pairs of models of
$T$).

\begin {definition}\label {foext}\em
We say that the extension property is first-order in $T$ iff for every formulas
$\phi(x,y),\psi(y,z)\in L$ the relation $Q_{\phi,\psi}$ defined by: $$Q_{\phi,\psi}(a)\mbox{\ iff}\
\phi(x,b)\mbox{ doesn't\ fork\ over}\ a\ \mbox{for\ every}\ b\models \psi(y,a)$$ is type-definable
(here $a$ can be an infinite tuple from $\CC$ whose sorts are fixed). We say that $T$ has \em wnfcp
\em if $T$ is low and the extension property is first-order in $T$.
\end {definition}

\begin{remark}\em
Recall that $T$ has the nfcp (non finite cover property) iff for every formula $\phi(x,y)\in L$
there exists $k<\omega$ such that every set $\{ \phi(x,a_i)\vert i\in I\}$ of instances of
$\phi(x,y)$ is consistent iff every subset of it of size $k$ is consistent. By a theorem of Shelah,
$T$ has nfcp iff $T$ is stable and $T^{eq}$ eliminates the quantifier $\exists^\infty$ [Sh, Chapter
2, Theorems 4.2, 4.4]. Moreover, if $T$ is stable then $T$ has the nfcp iff $T$ has the wnfcp
[BPV].
\end{remark}

\begin{fact}$[S1, Corollary\ 3.13]$\label{ext pcft}
Suppose the extension property is first-order in $T$. Then $T$ is PCFT.
\end{fact}

We say that an $A$-invariant set $\UU$ \em has finite $SU$-rank \em if $SU(a/A)<\omega$ for all
$a\in\UU$, and \em has bounded finite $SU$-rank \em if there exists $n<\omega$ such that
$SU(a/A)\leq n$ for all $a\in\UU$. The existence of a $\tau^f$-open set of bounded finite $SU$-rank
implies the existence of an $SU$-rank 1 formula (i.e. a weakly-minimal formula):

\begin{fact}\label{tau bounded SU}$[S0, Proposition\ 2.13]$
Let $\UU$ be an unbounded $\tau^f$-open set over some set $A$. Assume $\UU$ has bounded finite
$SU$-rank. Then there exist a set $B\supseteq A$ with $\vert B\backslash A\vert<\omega$ and
$\theta(x)\in L(B)$ of $SU$-rank 1 such that $\theta^\CC\subseteq \UU\cup acl(B)$.
\end{fact}

In [S1] the class of $\tilde\tau^f$-sets and its subclass of $\tilde\tau^f_{st}$-sets were
introduced. The class of $\tilde\tau^f$-sets is much wider than the class of basic $\tau^f$-open
sets. Here we look at the intermediate class of $\tilde \tau^f_{low}$-sets.

\begin{definition}\em
A relation $V(x,z_1,...z_l)$ is said to be a \em pre-$\tilde\tau^f$-set relation over $\emptyset$
\em if there are $\theta(\tilde x,x,z_1,z_2,...,z_l)\in L$ and $\phi_i(\tilde x,y_i)\in L$ for
$0\leq i\leq l$ such that for all $a,d_1,...,d_l$ from $\CC$ we have
$$V(a,d_1,...,d_l)\ \mbox{iff}\ \exists \tilde a\
 [\theta(\tilde a,a,d_1,d_2,...,d_l)\wedge\bigwedge^{l}_{i=0} (\phi_i(\tilde a,y_i)\ \mbox{forks\
over}\ d_1d_2...d_i)]\ $$ (for $i=0$ the sequence $d_1d_2...d_i$ is interpreted as $\emptyset$). If
each $\phi_i(\tilde x,y_i)$ is assumed to be low in $y_i$ , $V(x,z_1,...z_l)$ is said to be a \em
pre-$\tilde\tau_{low}^f$-set relation.
\end{definition}

\begin{definition}\em
1) A \em $\tilde\tau^f$-set over $\emptyset$ \em is a set of the form
$$\UU=\{a \vert\ \exists d_1,d_2,...d_l\
V(a,d_1,...,d_l)\}$$ for some pre-$\tilde\tau^f$-set relation $V(x,z_1,...z_l)$.

\noindent 2) A \em $\tilde\tau^f_{low}$-set over $\emptyset$ \em is a set of the form
$$\UU=\{a \vert\ \exists d_1,d_2,...d_l\
V(a,d_1,...,d_l)\}$$ for some pre-$\tilde\tau^f_{low}$-set relation $V(x,z_1,...z_l)$.
\end{definition}

\begin{remark}\label{tau_low_compact}\em
Every $\tilde\tau^f_{low}$-set is type-definable.
\end{remark}

\proof Let $\phi(x,y)\in L$ be low in $x$. Let $\Gamma_\phi(y,z)$ be the invariant relation defined
by $\Gamma_\phi(a,c)$ iff $\phi(x,a)$ divides over $c$. Then $\Gamma_\phi(y,z)$ is type-definable,
so the claim follows by compactness.

\section{The Theorem}
In this section $T$ is assumed to be a simple theory and we work in $\CC$ (so, $T$ not necessarily
eliminates imaginaries).

\begin{definition}\em
Let $\Theta=\{\theta_i(x_i,x)\}_{i\in I}$ be a set of $L$-formulas such that $\forall
 x\exists^{<\infty}x_i \theta_i(x_i,x)$ for all $i\in I$. Let $s$ be the sort of $x$. For $A\subseteq
\CC^s$, let $acl_\Theta(A)=\{b\vert\ \theta_i(b,a)\ \mbox{for\ some}\ \theta_i\in\Theta\
\mbox{and}\ a\in A\}$.
\end{definition}

\begin{definition}\em
An invariant set $\UU(x,y_1,...y_r)$ is said to be \em a generalized uniform family of
$\tilde\tau^f_{low}$-sets \em if there is a formula $\rho(\tilde
x,x,y_1,...,y_r,z_1,z_2,...,z_k)\in L$ and there are formulas $\psi_i(\tilde x,v_i),\mu_j(\tilde
x,w_j)\in L$ for $0\leq i\leq r$ and $1\leq j\leq k$ that are low in $v_i$ and low in $w_j$,
respectively, such that for all $a,d_1,...d_r$ we have $\UU(a,d_1,...d_r)$ iff $\exists \tilde
a\exists e_1...e_k$ $$\rho(\tilde a,a,d_1,...d_r,e_1,...e_k)\wedge [\bigwedge_{i=0}^r(\psi_i(\tilde
a,v_i)\ \mbox{forks\ over}\ d_1...d_i)]\wedge [\bigwedge_{j=1}^k(\mu_j(\tilde a,w_j)\ \mbox{forks\
over}\ d_1...d_re_1...e_j)].$$
\end{definition}

\begin{definition}\em
An invariant set $\FF(x,y_1,...y_r)$ is said to be \em a generalized uniform family of
$\tilde\tau^f_{low}$-closed sets \em if $\FF(x,y_1,...y_r)=\bigcap_i \neg\UU_i(x,y_1,...y_r)$,
where each $\UU_i(x,y_1,...y_r)$ is a generalized uniform family of $\tilde\tau^f_{low}$-sets.
\end{definition}

The following fact [S1, Theorem 8.7] is the key ingredient of our main theorem.

\begin{fact}\label {tilde top thm}
Assume the extension property is first-order in $T$. Let $\UU$ be an unbounded $\tilde\tau^f$-set
over $\emptyset$. Then there exists an unbounded $\tau^f$-open set $\UU^*$ over some finite set
$A^*$ such that $\UU^*\subseteq \UU$. In fact, if $V(x,z_1,...,z_l)$ is a pre-$\tilde\tau^f$-set
relation such that $\ \UU=\{a\vert \exists d_1...d_l V(a,d_1,...,d_l)\}$, and $\bar
d^*=(d^*_1,...,d^*_m)$ is any maximal sequence (with respect to extension) such that $\UU^*_{\bar
d^*}=\exists d_{m+1}...d_l V(\CC,d^*_1,...,d^*_m,d_{m+1},...,d_l)$ is unbounded, then $\UU^*_{\bar
d^*}$ is a $\tau^f$-open set over $d^*_1...d^*_m$.

\end{fact}

\begin{theorem}\label{main_thm}
Let $T$ be a countable simple theory in which the extension property is first-order. Assume:\\
1) $\Theta=\{\theta_i(x'_i,x)\}_{i<\omega}$ is a set of $L$-formulas such that $\forall
x\exists^{<\infty}x'_i
\theta_i(x'_i,x)$ for all $i<\omega$.\\
2) $\UU_0(x)$ is a non-empty $\tilde\tau^f_{low}$-set over $\emptyset$.\\
3) $\{F_n(x_n)\}_{n<\omega}$ is a family of $\emptyset$-invariant sets such that $F_n(\CC)\cap
acl(\emptyset)=\emptyset$ for all $n<\omega$.\\
4) For every $n<\omega$ and every variables $\bar y=y_1,...y_r$, let $\FF_n^{\bar y}(x_n,\bar y)$
be a generalized uniform family of $\tilde\tau^f_{low}$-closed sets such that
$F_n(\CC)\subseteq \FF_n^{\bar y}(\CC,\bar d)$ for all $\bar d$.\\
Now, assume that for all $a\in\UU_0$ there exists $b\in acl_{\Theta}(a)$ and $n<\omega$ such that
$b\in F_n(\CC)$. Then there is an unbounded $\tau^f_\infty$-open set $\UU^*$ over a finite tuple
$\bar d^*$ and variables $\bar y^*$ of the sort of $\bar d^*$, and $n^*<\omega$ such that
$$\UU^*\subseteq \FF_{n^*}^{\bar y^*}(\CC,\bar d^*)\cap acl_{\Theta}(\UU_0).$$
\end{theorem}

\proof First, we may assume $\Theta$ is closed downwards (i.e. if $\theta\in\Theta$ and
$\theta'\vdash \theta$ then $\theta'\in \Theta$; note that since $L$ is countable the closure of
$\Theta$ in this sense remains countable). Assume the conclusion of the theorem is false. To get a
contradiction, it will be sufficient to show the following.

\begin{subclaim}\label{main_subclaim}
For every non-empty $\tilde\tau^f_{low}$-set $\UU\subseteq \UU_0$ over $\emptyset$, every
$\theta\in \Theta$, and every $n<\omega$ there exists a non-empty $\tilde\tau^f_{low}$-set
$\UU^*\subseteq \UU$ over $\emptyset$ such that either $\neg\exists x'\theta(x',a)$ for all
$a\in\UU^*$ or for all $a\in\UU^*$ there exists $b\models \theta(x',a)$ with $b\not\in F_n(\CC)$.
\end{subclaim}

First, we show this is sufficient. Construct a decreasing sequence $(\UU_m \vert m<\omega)$ of
non-empty $\tilde\tau^f_{low}$-sets that begins at $\UU_0$, and for every $m<\omega$ the set
$\UU_{m+1}$ is obtained from $\UU_m$ by applying Subclaim \ref{main_subclaim} for an appropriate
pair $(\theta,n)$ (that corresponds to $m$ by a fixed bijection of $\Theta\times\omega$ with
$\omega$). By Remark \ref{tau_low_compact} and compactness $\bigcap \UU_m\neq \emptyset$, so there
exists $a^*\in \UU_0$ such that for all $\theta\in\Theta$ either $\neg\exists x'\theta(x',a^*)$ or
for every $n<\omega$ there exists $b_{n,\theta}\models\theta(x',a^*)$ such that
$b_{n,\theta}\not\in F_n(\CC)$. Now, by the assumption of the theorem there exist
$\theta(x',x)\in\Theta$, $b^*$ and $n^*<\omega$ such that $\theta(b^*,a^*)$ and $b^*\in
F_{n^*}(\CC)$. As $\Theta$ is closed downwards, there exists $\theta^*(x',x)\in\Theta$ such that
$\theta^*(x',x)\vdash \theta(x',x)$ and $\theta^*(x',a^*)$ isolates $tp(b^*/a^*)$ (as it is
algebraic). By the above property of $a^*$, there exists $b^{**}\models \theta^*(x',a^*)$ with
$b^{**}\not\in F_{n^*}(\CC)$; a contradiction to
the fact that $\theta^*(x',a^*)$ isolates $tp(b^*/a^*)$ and the assumption that $F_{n^*}(\CC)$ is $\emptyset$-invariant.\\


\noindent\textbf{Proof of Subclaim \ref{main_subclaim}} To show this, let $\UU$, $\theta$ and
$n<\omega$ be given. Let $V(x,z_1,...z_l)$ be a pre-$\tilde\tau^f_{low}$-set relation such that
$$\UU=\{a \vert\ \exists d_1,d_2,...d_l\ V(a,d_1,...,d_l)\}.$$ where $V$ is defined by:
$$V(a,d_1,...,d_l)\ \mbox{iff}\ \exists \tilde a\
 [\sigma(\tilde a,a,d_1,d_2,...,d_l)\wedge\bigwedge^{l}_{i=0} (\phi_i(\tilde a,t_i)\ \mbox{forks\
over}\ d_1d_2...d_i)]$$ for some $\sigma(\tilde x,x,z_1,z_2,...,z_l)\in L$ and $\phi_i(\tilde
x,t_i)\in L$ which are low in $t_i$ for $0\leq i\leq l$. Let $V_\theta$ be defined by: for all
$b,d_1,...,d_l\in \CC$, $$V_\theta(b,d_1,...,d_l)\ \mbox{iff}\ \exists a (\theta(b,a)\wedge
V(a,d_1,...,d_l)).$$ and let
$$\UU_\theta=\{b \vert\ \exists d_1,d_2,...d_l\ V_\theta(b,d_1,...,d_l)\}.$$ Since by the assumption $F_n(\CC)\cap
acl(\emptyset)=\emptyset$, we may assume $\UU_\theta\cap acl(\emptyset)=\emptyset$ and $\UU_\theta$
is non-empty. Now, let $\bar d^*=(d_1^*,...,d_m^*)$ be a maximal sequence, with respect to
extension ($0\leq m\leq l$), such that $$\tilde V_\theta(x')\equiv \exists d_{m+1},d_{m+2},...d_l
V_\theta(x',d^*_1,...d^*_m,d_{m+1},...d_l)$$ is non-algebraic. We may assume $m<l$ (by choosing $V$
appropriately). By Fact \ref{tilde top thm}, $\tilde V_\theta(\CC)$ is an unbounded basic
$\tau^f_\infty$-open set over $\bar d^*$. Since we assume the conclusion of the theorem is false,
$\tilde V_\theta(\CC)\not\subseteq \FF_{n}^{\bar y^*}(\CC,\bar d^*)$ where $\bar
y^*=y_1^*,...,y_m^*$ has the same sort as $\bar d^*$. Now, let each $\UU_{s,n}(x_n,\bar y^*)$ for
$s<\alpha$ be a generalized uniform family of $\tilde \tau^f_{low}$-sets such that $\FF_n(x_n,\bar
y^*)=\bigcap_{s<\alpha} \neg\UU_{s,n}(x_n,\bar y^*)$. Let $b^*\in\tilde
V_\theta(\CC)\backslash\FF_{n}^{\bar y^*}(\CC,\bar d^*)$. So, there exists $s^*<\alpha$ such that
$b^*\in\UU_{s^*,n}(\CC,\bar d^*)$. Let $\rho(\tilde x',x_n,y^*_1,...,y^*_m,z'_1,z'_2,...,z'_k)\in
L$ and let $\psi_i(\tilde x',v_i),\mu_j(\tilde x',w_j)\in L$ for $0\leq i\leq m$ and $1\leq j\leq
k$ be low in $v_i$ and low in $w_j$ respectively, such that for all $b,d_1,...d_m$ we have
$\UU_{s^*,n}(b,d_1,...d_m)$ iff $\exists \tilde b\exists e_1...e_k$ $$\rho(\tilde
b,b,d_1,...d_m,e_1,...e_k)\wedge [\bigwedge_{i=0}^m(\psi_i(\tilde b,v_i)\ \mbox{forks\ over}\
d_1...d_i)]\wedge [\bigwedge_{j=1}^k(\mu_j(\tilde b,w_j)\ \mbox{forks\ over}\
d_1...d_me_1...e_j)].$$ Now, let $d^*_{m+1},...d^*_l$ and $a^*,\tilde a^*$ and
$E^*=(e^*_1,...,e^*_k)$ and $\tilde b^*$ be such that
$$\theta(b^*,a^*)\wedge\sigma(\tilde a^*,a^*,d^*_1,d^*_2,...,d^*_l)\wedge\bigwedge^{l}_{i=0}
 (\phi_i(\tilde a^*,y_i)\ \mbox{forks\ over}\ d^*_1d^*_2...d^*_i)\ \ (*1)$$ and $$\rho(\tilde
 b^*,b^*,d^*_1,..d^*_m,e^*_1,..e^*_k)\ (*2)$$ and $$[\bigwedge_{i=0}^m(\psi_i(\tilde
b^*,v_i)\ \mbox{forks\ over}\ d^*_1...d^*_i)]\wedge [\bigwedge_{j=1}^k(\mu_j(\tilde b^*,w_j)\
\mbox{forks\ over}\ d^*_1...d^*_me^*_1...e^*_j)]\ (*3).$$ By maximality of $\bar d^*$, we know
$b^*\in acl(\bar d^*d^*_{m+1})$. Thus, by taking a non-forking extension of $tp(\tilde
b^*E^*/acl(\bar d^*d^*_{m+1}))$ over $acl(d^*_1...d^*_la^*\tilde a^*)$ we may assume $E^*$ is
independent from $d^*_1...d^*_la^*\tilde a^*$ over $\bar d^*d^*_{m+1}$ and $(*1)$, $(*2)$ and
$(*3)$ still hold. We conclude that $$\bigwedge^{l}_{i=m+1} (\phi_i(\tilde a^*,t_i)\ \mbox{forks\
over}\ d^*_1d^*_2...d^*_iE^*).$$ Now, we define the $\tilde\tau^f_{low}$-set $\UU^*$. First, define
a relation $V^*$ by:
$$V^*(a,d_1,...d_m,e_1,...e_k,d_{m+1},..d_l)\ \mbox{iff}\ \exists \tilde a, b,\tilde b(\theta^*\wedge
V^*_0\wedge V^*_1\wedge V^*_2),$$ where $\theta^*$ is defined by: $\theta^*(\tilde a,b,\tilde
b,a,d_1,..d_m,e_1,...e_k,d_{m+1},..d_l)$ iff\
$$\theta(b,a)\wedge\sigma(\tilde a,a,d_1,d_2,...,d_l)\wedge \rho(\tilde b,b,d_1,...d_m,e_1,...,e_k),$$
$V^*_0$ is defined by: $V^*_0(\tilde a,\tilde b,d_1,...d_m)$ iff
$$\bigwedge^{m}_{i=0}(\phi_i(\tilde a,t_i)\vee \psi_i(\tilde b,v_i)\ \mbox{forks\ over}\ d_1d_2...d_i),$$ \noindent
$V^*_1$ is defined by $V_1(\tilde b,d_1,..d_m,e_1,...e_k)$ iff
$$\bigwedge_{j=1}^k(\mu_j(\tilde b,w_j)\ \mbox{forks\ over}\ d_1...d_me_1...e_j),\ \mbox{and}$$
$V^*_2$ is defined by $V_2(\tilde a,d_1,...d_m,e_1,...e_k,d_{m+1},..d_l)$ iff
$$\bigwedge^{l}_{i=m+1} (\phi_i(\tilde a,t_i)\ \mbox{forks\ over}\ d_1d_2...d_ie_1...e_k).$$
Note that $V^*$ is a pre-$\tilde\tau^f_{low}$-set. Let $$\UU^*=\{a \vert \exists
d_1,..d_m,e_1,...e_k,d_{m+1},..d_l\ V^*(a,d_1,..d_m,e_1,...e_k,d_{m+1},...d_l)\}.$$ By the
definition of $\UU^*$, $\UU^*\subseteq \UU$. $\UU^*$ is a $\tilde\tau^f_{low}$-set using Remark
\ref{low_disjunction}. By the construction, $\UU^*\neq\emptyset$. Now, let $a\in\UU^*$. By the
definition of $\UU^*$, there are $\tilde b,b,d_1,...d_m,e_1,...e_k$ such that $\theta(b,a)$,
$\rho(\tilde b,b,d_1,...d_m,e_1,...,e_k)$, $$\bigwedge_{i=0}^m(\psi_i(\tilde b,v_i)\ \mbox{forks\
over}\ d_1,...d_i)\ \mbox{and}\ \bigwedge_{j=1}^k(\mu_j(\tilde b,w_j)\ \mbox{forks\ over}\
d_1,...d_me_1...e_j).$$ Thus $\UU_{s^*,n}(b,d_1...d_m)$ and therefore $\neg\FF_{n}^{\bar
y^*}(b,d_1...d_m)$. Hence $b\not\in F_n$ as required.

\section{Applications}

In this section we show some applications of Theorem \ref{main_thm}. In fact, we will show several
instances of this theorem that are apparently new even for stable theories. In this section $T$ is
assumed to be a simple theory and we work in $\CC$.\\

We start by pointing out that theorem \ref{main_thm} generalizes [S1, Theorem 9.4] that is one of
the essential steps towards the proof of supersimplicity of countable simple unidimensional
theories with elimination of hyperimaginaries. First recall the following definitions from [S1] of
stable-independence and the $SU_{se}$-rank.

\begin{definition}\em
For $a\in \CC$, $A,B\subseteq \CC$,  $\sfork{a}{B}{A}$ if for some stable $\phi(x,y)\in L$, there
is $b$ in $A\cup B$ and $a'\in\phi(\CC,b)\cap dcl(Aa)$ such that $\phi(x,b)$ forks over $A$.
\end{definition}

\begin{definition}\em
The $SU_{se}$-rank of $tp(a/A)$ is defined by induction on $\alpha$: if $\alpha=\beta+1$,
$SU_{se}(a/A)\geq \alpha$ if there exist $B_1\supseteq B_0\supseteq A$ such that
$\sfork{a}{B_1}{B_0}$ and $SU_{se}(a/B_1)\geq\beta$. For limit $\alpha$, $SU_{se}(a/A)\geq\alpha$
if $SU_{se}(a/A)\geq\beta$ for all $\beta<\alpha$.
\end{definition}

\begin{remark}\label{symmetry_s}\em
In [S1, Lemma 6.8] it is proved that in a simple theory, in which $Lstp=stp$ over sets,
$\snonforkempty{}{}$ is symmetric. In fact, $\snonforkempty{}{}$ is symmetric in any simple theory.
Thus for any simple theory, if $s_0$ and $s_1$ are finite tuples of sorts and $n<\omega$ then the
set $\FF^{s_0,s_1}_n$ defined by $$\FF^{s_0,s_1}_n=\{ (a,A)\in \CC^{s_0}\times \CC^{s_1} \vert\
SU_{se}(a/A)<n\}$$ is a generalized uniform family of $\tilde\tau^f_{low}$-closed sets.
\end{remark}

\proof To prove that $\snonforkempty{}{}$ is symmetric, first recall [S1, Claim 6.5]:

\begin{fact}\label{stable claim}
Let $T$ be simple. Let $\phi(x,y)\in L$ be stable. Assume $\nonfork{a}{b}{A}$ and
$\nonfork{a'}{b}{A}$ and $Lstp(a/A)=Lstp(a'/A)$. Then $\phi(a,b)$ iff $\phi(a',b)$.
\end{fact}

By the proof of symmetry of stable-independence [S1, Lemma 6.8] it will be sufficient to prove Fact
\ref{stable claim} with the weaker assumption $stp(a)=stp(a')$ instead of the assumption
$Lstp(a)=Lstp(a')$ (we may clearly assume $A=\emptyset$). Indeed to prove this assume
$stp(a)=stp(a')$. Now, for every complete type $q\in S(\emptyset)$ let $E_q$ be the equivalence
relation defined by: $E_q(a,a')$ iff "for every $b\models q$ that is independent from $aa'$ we have
[$\phi(a,b)$ iff $\phi(a',b)$]".Then $E_q$ Stone-open. By Fact \ref{stable claim}, equality of the
Lascar strong type refines $E_q$. Thus $E_q$ is a $\emptyset$-definable finite equivalence relation
(as a bounded Stone-open equivalence relation is definable [S4, Lemma 7]). Now, by the assumption
that $stp(a)=stp(a')$, $E_q(a,a')$ for all complete $q$. Thus, by extension we get that for every
$b$, if each of $a$ and $a'$ is independent from $b$, then $\phi(a,b)$ iff $\phi(a',b)$.

We explain now the last phrase. We need to show that $\neg\FF^{s_0,s_1}_n$ is a disjunction of
invariant sets, each of which is a generalized uniform family of $\tilde\tau^f_{low}$-sets for all
$s_0,s_1$ and $n$ as above. Indeed, by symmetry of $\sforkempty{}{}$, $\neg\FF^{s_0,s_1}_n(a,A)$
iff there are $b_1,c_1,...,b_n,c_n$ such that $\sfork{c_i}{a}{Ab_1c_1...b_{i-1}c_{i-1}b_i}$ for all
$1\leq i\leq n$. By the definition of $\sforkempty{}{}$, this can easily seen to be equivalent to
a disjunction of the required form (since any stable $\phi(x,y)\in L$ is low both in $x$ and in $y$).\\

\noindent For an $A$-invariant set $V$, let $acl_1(V)=\{a'\vert\ a'\in acl(a)\ \mbox{for\ some}\
a\in V^1\}$. The following corollary generalizes [S1, Theorem 9.4].

\begin{corollary}\label{cor1}
Let $T$ be a countable simple theory in which the extension property is first-order. Let $\UU_0$ be
a non-empty $\tilde\tau^f_{low}$-set. Assume for every $a\in\UU_0$ there exists $a'\in
acl(a)\backslash acl(\emptyset)$ such that $SU_{se}(a')<\omega$. Then there exists an unbounded
$\tau_{\infty}^f$-open set $\UU\subseteq acl_1(\UU_0)$ over a finite set such that $\UU$ has
bounded finite $SU_{se}$-rank.
\end{corollary}

\proof Let $x$ be the variable of $\UU_0$, so $\UU_0=\UU_0(x)$. Let $$\Theta=\{\theta(x',x)\vert\
\exists^{<\infty} x' \theta(x',x), x'\ \mbox{any variable}\}.$$ Let $\SS$ be the set of sorts. Let
$I:\omega\rightarrow \SS\times \omega$ be a bijection, $I_1,I_2$ the projections of $I$ to the
first and second coordinate, respectively. Now, for each $n<\omega$ let $F_n=\{a\in
\CC^{I_1(n)}\backslash acl(\emptyset)\ \vert SU_{se}(a)<I_2(n)\}$. Now, for every finite tuple of
variables $Y$ and $n<\omega$ let $s(Y)$ be the finite sequence of sorts of $Y$ and let
$$\FF^Y_n=\{(a,A)\in \CC^{I_1(n)}\times \CC^{s(Y)}\vert\ SU_{se}(a/A)<I_2(n)\}.$$ Now, by the
definition of the $SU_{se}$-rank, $F_n(\CC)\subseteq\FF^Y_n(\CC,A)$ for every $n<\omega$ and every
$Y,A$. By Remark \ref{symmetry_s}, $\FF^Y_n$ is a generalized uniform family of
$\tilde\tau^f_{low}$-closed sets for all $Y,n$. By our assumptions, we see that the assumptions of
Theorem \ref{main_thm} hold for $\UU_0(x)$, $\Theta$ ,$\{F_n\}_n$ and $\{\FF^Y_n\}_{Y,n}$ and thus
by its conclusion we are done.

\begin{corollary}\label {cor2}
Let $T$ be a countable theory with wnfcp. Let $\UU_0$ be an unbounded $\tilde\tau^f$-set over
$\emptyset$ of finite $SU$-rank. Then there exists a finite set $A$ and an $SU$-rank 1 formula
$\theta\in L(A)$ such that $\theta^\CC\subseteq \UU_0\cup acl(A)$.
\end{corollary}

\proof First, by modifying $\UU_0$, we may assume $\UU_0\cap acl(\emptyset)=\emptyset$. Let
$\Theta=\{x'=x\}$, $\UU_0(x)=\UU_0$. Let $s(x)$ be the sort of $x$. Now, for each $n<\omega$ let
$$F_n=\{a\in \CC^{s(x)}\backslash acl(\emptyset)\ \vert SU(a)<n\}.$$ For every finite tuple
of variables $Y$ and $n<\omega$ let $s(Y)$ be the finite sequence of sorts of $Y$ and let
$$\FF^Y_n=\{(a,A)\in \CC^{s(x)}\times \CC^{s(Y)}\vert\ SU(a/A)<n\}.$$ By symmetry of
forking and the assumption that $T$ is low, each $\FF^Y_n$ is a generalized uniform family of
$\tilde\tau^f_{low}$-closed sets. Clearly, $F_n(\CC)\subseteq \FF_n^Y(\CC,A)$ for every $n<\omega$
and every $Y,A$. By our assumption, the assumptions of Theorem \ref{main_thm} are satisfied for
$\UU_0$, $\Theta$, $\{F_n\}_n$ and $\{\FF^Y_n\}_{Y,n}$ and thus by its conclusion there exists an
unbounded $\tau^f_{\infty}$-open set $\UU^*\subseteq \UU_0$ over a finite set $A_0$ and $\UU^*$ has
bounded finite $SU$-rank. By Fact \ref{tau bounded SU}, there exists a finite set $A\supseteq A_0$
and there exists a $SU$-rank 1 formula $\theta\in L(A)$ such that $\theta^\CC\subseteq \UU^*\cup
acl(A)$.

\begin{corollary}
Let $T$ be a countable theory with wnfcp. Let $\UU_0$ be a non-empty $\tilde\tau^f$-set over
$\emptyset$. Assume for every $a\in\UU_0$ there exists $a'\in acl(a)\backslash acl(\emptyset)$ such
that $SU(a')<\omega$. Then there exists a finite set $A$ and an $SU$-rank 1 formula $\theta\in
L(A)$ such that $\theta^\CC\subseteq acl_1(\UU_0)\cup acl(A)$.
\end{corollary}

\proof Just like the proof of Corollary \ref{cor2}.

\section{Dichotomies for countable theories with the wnfcp} In this section we show that the dichotomy
[S1, Theorem 5.5] implies a strong dichotomy between essential 1-basedness and supersimplicity in
the case $T$ is a countable wnfcp theory that eliminates hyperimaginaries. Before we state the
above dichotomy for the special case of the $\tau^f$-topologies (simplified version), let us recall
the basic definitions. In this section $T$ is assumed to be simple and we work in $\CC=\CC^{eq}$.\\

First, let us fix some notations and terminology. Let $V,W$ be invariant sets. We say that $V$ is
generated over $W$ by a small set $B$ if $V\subseteq dcl(W\cup B)$. We say that $V$ is generated
over $W$ if it is generated over $W$ by some small set. If $V$ is $A$-invariant, we say that $V$ is
(almost) $W$-internal over $A$ if for every $a\in V$ there exists $B\supseteq A$, over which $W$ is
invariant, that is independent from $a$ over $A$ and there exists a tuple $\bar c$ of realizations
of $W$ such that $a\in dcl(B,\bar c)$ ($a\in acl(B,\bar c)$, respectively). If we say that $V$ is
$W$-internal (without specifying over what set) then we mean that $V$ is $W$-internal over the set
that $V$ comes with (e.g. in case it is a partial type, we consider it with its specified
parameters). Note that if both $V$ and $W$ are $A$-invariant then for all $B,C\supseteq A$, $V$ is
(almost) $W$-internal over $B$ iff $V$ is (respectively, almost) $W$-internal over $C$.


\begin{definition}\label {def ess-1-based}\em
A type $p\in S(A)$ is said to be \em essentially 1-based by means of the $\tau^f$-topologies \em if
for every finite tuple $\bar c$ from $p$ and for every type-definable $\tau^f$-open set $\UU$ over
$A\bar c$, the set $\{a\in \UU \vert\ Cb(a/A\bar c)\not\in bdd(aA)\}$ is nowhere dense in the
Stone-topology of $\UU$.
\end{definition}

We state now [S1, Theorem 5.5] for the $\tau^f$-topologies (in fact, it is a special case of it
when working over constants). Also, as indicated in the end of the proof of this fact, the
finite-$SU$-rank $\tau^f$-open set we obtained is almost $p_0$-internal.

\begin{fact}\label{dichotomy thm}
Let $T$ be a countable simple theory with PCFT that eliminates hyperimaginaries. Let $p_0$ be a
partial type over $\emptyset$ of $SU$-rank 1. Then, either there exists an unbounded $\tau^f$-open
set over some countable set that is almost internal to $p_0$ (in particular, has finite-$SU$-rank)
or every type $p\in S(A)$, with $A$ countable, that is internal in $p_0$ is essentially 1-based by
means of the $\tau^f$-topologies.
\end{fact}

\begin{theorem}\label{dichotomy thm_wnfcp}
Let $T$ be a countable theory with wnfcp that eliminates hyperimaginaries. Let $p$ be a partial
type over $\emptyset$ of $SU$-rank 1 . Then, either

\noindent 1) every type $q\in S(A)$, with $A$ countable, that is internal in $p$ is essentially
1-based by means of the $\tau^f$-topologies, or

\noindent 2) there exists a weakly-minimal definable set (in $L(\CC)$) that is generated over
$p(\CC)$.
\end{theorem}

\proof Assume 1) is false. By Fact \ref{dichotomy thm}, there exists an unbounded type-definable
$\tau^f$-open set $\UU$ over some countable set $A$ such that $tp(a/A)$ is almost $p$-internal for
every $a\in \UU$.

\begin{subclaim}
There exists an unbounded type-definable $\tau^f$-open set $\UU^*$ over $A$ that is generated over
$p(\CC)$.
\end{subclaim}

\proof By [WB] or [S2, Corollary 4.9], for every $a\in \UU\backslash acl(A)$ there exists $a'\in
dcl(aA)\backslash acl(A)$ such that $tp(a'/A)$ has fundamental system of solutions over $p(\CC)$
(i.e. $tp(a'/A)$ is generated over $p(\CC)$ by a set of realizations of $tp(a'/A)$ together with
$A$.) In particular, there exists a (finite) set $A'$ of realizations of $tp(a'/A)$ that is
independent from $a'$ over $A$ and tuple $\bar c$ of realizations of $p$ such that $a'\in
dcl(A'A\bar c)$. For every
$A$-definable functions $f,g$ let\\

\noindent $F_{f,g}=\{a\in \UU \vert f(a)=g(\bar b,\bar c)\not\in acl(A)\ \mbox{for\ some\ } \bar
b,\bar c\ \mbox {with} \nonfork{f(a)}{\bar b}{A},$

\noindent $\ \ \ \ \ \ \ \ \ \ \ \ \ \ \ \ \mbox{where\ } \bar c\ \mbox{is\ a tuple\ of\
realizations\ of\ } p,\ \mbox{and}\ \bar b\ \mbox{is\ a\ tuple\ of\ realizations\ of}\ tp(f(a)/A)\}$.\\

\noindent By Remark \ref{exists generic}(3), each $F_{f,g}$ is $\tau^f$-closed over $A$. Thus, by
Baire category theorem for the $\tau^f$-topology (by Remark \ref{exists generic}(2),
($\UU\backslash acl(A), \tau^f)$ is a Baire space) there are $A$-definable functions $f^*,g^*$ such
that $F_{f^*,g^*}$ has non-empty interior in the $\tau^f$-topology over $A$. By Fact \ref{ext pcft}
there exists an unbounded type-definable $\tau^f$-open set $\UU^*$ over $A$ such that for every
$a\in\UU^*$ there exists a tuple $\bar b$ of realizations of $tp(a/A)$ that is independent from $a$
over $A$ such that $a=g^*(\bar b,\bar c)$ for some tuple $\bar c$ of realizations of $p$. The
subclaim follows now directly from [S2, Theorem 3.7]:

\begin {fact}\label {gen}
Let $p\in S(\emptyset)$ and let $\RR$ be $\emptyset$-invariant. Suppose the internality of $p$ in
$\RR$ is witnessed by a generic parameter whose type $q$ is almost-$\RR$-internal. Then $p$ is
generated over $\RR$ by a set of realizations of $q$.
\end {fact}


Now, as $\UU^*$ has bounded finite $SU$-rank (the bound is determined by $g^*$), by Fact \ref{tau
bounded SU}, there exists an $SU$-rank 1 formula $\theta(x,b)$ such that $\theta(\CC,b)\subseteq
\UU^*\cup acl(Ab)$. Thus 2) follows.

\subsection{A trichotomy for countable theories with the nfcp}

Here we prove a trichotomy for countable theories with the nfcp. In this subsection we work in a
large saturated model $\CC=\CC^{eq}$ of a simple theory $T$ with
elimination of hyperimaginaries unless stated otherwise.\\

We begin with some standard terminology and remarks. For a definable set $D$ over $A$ we denote by
$D^*$ the induced structure on $D$ over $A$, namely, $D^*$ is the set $D$ equipped with all
$A$-definable relations in $\CC$ that are subsets of $D^n$ for some $n$. Then, easily $D^*$ has
elimination of quantifiers and therefore saturated.

\begin{definition}\em
Let $D$ be a type-definable set over a set $A$. We say that $D$ is 1-based if for every finite
tuple $\bar a$ of realizations of $D$ and and set $B\supseteq A$, we have $Cb(\bar a/B)\in acl(\bar
aA)$. A type-definable group $G$ over $A$ is said to be 1-based if its underlying set is.
\end{definition}

\begin{remark}\em\label{dichotomy_remark}
1) A type-definable set $D$ over $A$ is 1-based iff $\bar a$ is independent from $\bar a'$ over
$acl(A\bar a)\cap acl(A\bar a')$ for every finite tuples $\bar a$ and $\bar a'$
from $D$.\\
2) Let $D$ be a definable set over $A$. Then

\noindent i) if $T$ is stable (simple), so is $Th(D^*)$.

\noindent ii) if $D^*$ is 1-based then $D$ is 1-based (as a type-definable set).

\noindent iii) if $D$ is stably-embedded (e.g. $T$ is stable), and $p$ is a partial type of $D^*$
then $RM_{D^*}(p)=RM(p_D)$ (where $p_D$ is just the conjunction of $p$ with appropriate power of
$D$, $RM$ is the usual Morley rank in $\CC$, and $RM_{D^*}$ is the Morley rank in $D^*$).
\end{remark}

\begin{lemma}\label{morley_rank}
Assume $L$ is countable and  $\theta(\CC)\subseteq acl(p(\CC))$, where $p$ is any partial type over
$\emptyset$ and $\theta(x)\in L$ is non-algebraic. Then

\noindent 1) there exists a $\emptyset$-definable $\theta^*(x)\vdash \theta(x)$ and
$\emptyset$-definable functions $f,g$ and $n<\omega$ such that $f[\theta^*(\CC)\backslash
acl(\emptyset)]\subseteq g[{p^n}(\CC)]$ and $f[\theta^*(\CC)]$ is non-algebraic, and

\noindent 2) if $p$ is minimal then $f[\theta^*(\CC)]$ has ordinal Morley rank and thus contains a
strongly-minimal formula.
\end{lemma}

\proof For every $a\in \theta(\CC)\backslash acl(\emptyset)$ there exist $n<\omega$ and $\bar c\in
p^n(\CC)$ such that $a\in acl(\bar c)$. Let $e=Cb(\bar c/a)$. Now, by elimination of
hyperimaginaries there exists $e^*\in acl(a)\cap dcl(p(\CC))\backslash acl(\emptyset)$. Let
$e^{**}=\{e' \vert tp(e'/a)=tp(e^*/a)\}$  ($e^{**}$ is an imaginary element). Then, clearly
$e^{**}\in dcl(a)\cap dcl(p(\CC))\backslash acl(\emptyset)$. For any appropriate
$\emptyset$-definable functions $f,g$ let $$F_{f,g}=\{a\in \theta(\CC) \vert\ \exists\bar
c\subseteq p(\CC)\ [f(a)=g(\bar c)\not\in acl(\emptyset)] \}.$$ So, $\{F_{f,g}\}_{f,g}$ is a
countable family of Stone-closed sets that covers $\theta(\CC)\backslash acl(\emptyset)$ and thus
by Baire category theorem for the Stone-topology of $\theta(\CC)\backslash acl(\emptyset)$ we get
the required formula $\theta^*\in L$ and $\emptyset$-definable functions $f,g$ as in 1). To prove
2), assume that $p$ is minimal. Then, by induction on $n$, we easily get that for every countable
set $A$ the number of (complete) types of realizations of $p^n$ over $A$ is countable. Thus by 1),
for every countable set $A$ the number of complete types over $A$ extending $f[\theta^*(\CC)]$ is
countable. Therefore $f[\theta^*(\CC)]$ has
ordinal Morley rank.\\

We will be using the following two important facts. The first one is Buechler's dichotomy for
minimal types (see [P1, Corollary 3.3]).

\begin{fact}\label{buechler_dichotomy}
Let $T$ be superstable and let $p\in S(A)$ be a minimal type. Then either $p$ is 1-based or
$RM(p)=1$.
\end{fact}

The second fact is Wagner's result [W] on analysis in 1-based types in simple theories (it
generalizes previous results of Hrushovski and Chazidakis).

\begin{fact}\label{wagner_fact}
Let $T$ be any simple theory and work with hyperimaginaries. Assume $p\in S(A)$ is analyzable in an
$A$-invariant family of 1-based types. Then $p$ is 1-based.
\end{fact}

\begin{theorem}\label{trichotomy thm_nfcp}
Let $T$ be a countable theory with nfcp. Let $p\in S(\emptyset)$ be minimal. Then, either

\noindent 1) every type $q\in S(A)$, with $A$ countable, that is internal in $p$ is essentially
1-based by means of the $\tau^f$-topologies, or

\noindent 2) there is an infinite definable 1-based group of finite $D$-rank that is $p$-internal,
or

\noindent 3) there exists a strongly minimal definable set that is $p$-internal.

\end{theorem}

\proof Assume 1) is false. By Theorem \ref{dichotomy thm_wnfcp}, there exists a weakly-minimal
formula $\theta(x,b)$ that is $p$-generated and in particular $p$-internal (in the stable case an
invariant set is $p$-internal iff it is $p$-generated). First, assume $\theta(\CC,b)\subseteq
acl(p(\CC)\cup b)$. Then by Lemma \ref{morley_rank}, there exists a strongly-minimal formula
$\phi\in L(\CC)$ that is $p$-internal (even generated over $p^\CC$). Thus, we may assume
$\theta(\CC,b)\not\subseteq acl(p^\CC\cup b)$. Let $a\in \theta(\CC,b)\backslash acl(p^\CC\cup b)$.
Let $q=tp(a/acl(b))$ and let $\Gamma=Aut(q^\CC/p^\CC\cup acl(b))$. We will be using the following
fact [S2, Theorem 2.9], with its proof, which for simplicity we state for a special case. In the
following, for a set $S$, possibly large, we let $DCL(S)$ be the set of all elements in $\CC$ that
are fixed by any automorphism that fixes $S$ pointwise; we say that a set $V$ is controlled by $B$
over $S$, if $V\subseteq DCL(B\cup S)$.

\begin {fact}\label{bind_grp}
Let $T$ be any simple theory. Let $Q$ be a stably-embedded type-definable set over $\emptyset$ and
let $q\in S(\emptyset)$. Suppose there exists a set $B\subseteq DCL(q^\CC\cup Q)$ with $tp(B)\vdash
Lstp(B)$ such that $q^\CC$ is controlled by $B$ over $Q$. Then $\Gamma=Aut(q^\CC/Q)$ is
type-definable with its action on $q^\CC$ over $\emptyset$.
\end {fact}

\begin{remark}\label{bind_remark} \em
It is well known that in a stable theory if $q$ is $Q$-internal then there is always a set of
realizations $B$ of $q$ such that $q(\CC)\subseteq dcl(Q,B)$, in particular,  $q$ is controlled by
$B$ over $Q$; if $q$ is stationary then $B$ can be taken to be a finite initial segment of a Morley
sequence of $q$ and clearly $tp(B)\vdash Lstp(B)$. Now, $\Gamma$ in Fact \ref {bind_grp} can be
interpreted in the following way. As $Q$ is a type-definable stably-embedded set, there exists a
partial type $\sum_Q(Y,Y')$ expressing that $Y,Y'$ are $Q$-conjugate, for $Y,Y'\models tp(B)$. Now,
let $\Gamma_{B^2/Q}(Y,Y')$ be the type expressing that $tp(Y)=tp(Y')=tp(B)$ and $\sum_Q(Y,Y')$.
Now, by definition, $\sigma\in \Gamma=Aut(q^\CC/Q)$ iff $\sigma$ is the restriction to $q^\CC$ of
some automorphism of $\CC$ that fixes $Q$ pointwise. As $q$ is controlled by $B\subseteq
DCL(q^\CC\cup Q)$ over $Q$, it is not hard to show (see [S2, Theorem 2.9] proof) that $\Gamma$ can
be interpreted as $\Gamma_{B^2/Q}/E$ for certain $\emptyset$-definable equivalence relation $E$.
\end{remark}

By Remark \ref{bind_remark} and the fact that $q(x)\vdash \theta(x,b)$, there is an infinite
type-definable group $G$ over $acl(b)$ that is isomorphic to $\Gamma$ such that for some
$acl(b)$-definable equivalence relation $E$ and some $n<\omega$, we have $G\subseteq \theta(\CC,
b)^n/E$. Now, by stability of $T$, $G$ is an intersection of definable groups over $acl(b)$ [H1,
Theorem 2]. By compactness, there is an infinite $acl(b)$-definable group $G_0$ that is
$p$-internal and has finite $D$-rank. By Fact \ref {buechler_dichotomy} and Remark
\ref{dichotomy_remark}, 2)i) applied to the induced structure $G_0^*$ on $G_0$ over $acl(b)$, every
minimal type $r$ in $G^*_0$ is either 1-based or of Morley rank 1. Thus if 3) fails, then any such
$r$ is 1-based in $G^*_0$ by Remark \ref{dichotomy_remark}, 2)iii) and stability of $T$. As $G^*_0$
has finite $SU$-rank, we conclude, when working in $G^*_0$, that every non-algebraic type is
non-orthogonal to a minimal type, and therefore any type in $G^*_0$ is analyzable in 1-based types.
By Fact \ref{wagner_fact}, $G^*_0$ is 1-based. By Remark \ref{dichotomy_remark}, 2)ii), $G_0$ is
1-based.

\end{document}